\documentclass[10pt]{amsart}
\usepackage{amsbsy,amssymb,amscd,amsfonts,latexsym,amstext,delarray,
amsmath,amsthm,graphicx}

\newtheorem{The}{Theorem}[section]

\newtheorem{Pro}[The]{Proposition}

\newtheorem{Def}[The]{Definition}

\numberwithin{equation}{section}

\def\proof{\vspace{2ex}\noindent{\bf Proof.} }
\def\endproof{\relax\ifmmode\expandafter\endproofmath\else
\unskip\nobreak\hfil\penalty50\hskip.75em\hbox{}\nobreak\hfil\bull
{\parfillskip= 0pt \finalhyphendemerits= 0 \bigbreak}\fi}
\def\endproofmath$${\eqno\bull$$\bigbreak}
\def\bull{\vbox{\hrule\hbox{\vrule\kern3pt\vbox{\kern6pt}\kern3pt\vrule}
\hrule}}

\newcommand{\ba}{\begin{eqnarray}}
\newcommand{\na}{\end{eqnarray}}
 \def\t{\theta}
\def\fl{\forall}

\newcommand{\scr}{\mathcal}

 \def\Hc{{\mathcal H}}
 \def\g{\gamma}

\def\ve{\varepsilon}

\def\Qc{{\mathcal Q}}

 \def\Ec{{\mathcal E}}

 \def\a{\alpha}
 \def\b{\beta}

 \def\g{\gamma}

 \def\t{\theta}
 \def\ve{\varepsilon}

\newcommand{\C}{\mathbb{C}}
\newcommand{\bG}{\mathbb{G}}

\newcommand{\N}{\mathbb{N}}
\renewcommand{\P}{\mathbb{P}}
\newcommand{\Q}{\mathbb{Q}}
\newcommand{\R}{\mathbb{R}}

\newcommand{\Z}{\mathbb{Z}}

\newcommand{\fg}{{\mathfrak{g}}}

\newcommand{\fO}{{\mathcal{O}}}

\newcommand{\cF}{\scr{F}}

\newcommand{\cM}{\scr{M}}

\newcommand{\cU}{\scr{U}}

\newcommand{\ie}{{\it i.e.\/}\ }

\newcommand{\cf}{{\it cf.\/}\ }

\begin{document}

\title{Renormalization and motivic Galois theory}

\author[Connes]{Alain Connes}
\author[Marcolli]{Matilde Marcolli}
\address{A.~Connes: Coll\`ege de France, IHES, and Vanderbilt University}
\email{alain\@@connes.org}
\address{M.~Marcolli: Max--Planck Institut f\"ur Mathematik
Bonn, Germany} \email{marcolli\@@mpim-bonn.mpg.de}

\maketitle

\begin{abstract}
We investigate the nature of divergences in quantum field theory,
showing that they are organized in the structure of a certain
{\em`` motivic Galois group''} $U^*$, which is uniquely determined
and universal with respect to the set of physical theories. The
renormalization group can be identified canonically with a one
parameter subgroup of $U^*$. The group $U^*$ arises through a
Riemann--Hilbert correspondence. Its representations classify
equisingular flat vector bundles, where the equisingularity
condition is a geometric formulation of the fact that in quantum
field theory the counterterms are independent of the choice of a
unit of mass. As an algebraic group scheme, $U^*$ is a semi-direct
product by the multiplicative group $\bG_m$ of a pro-unipotent
group scheme whose Lie algebra is freely generated by one
generator in each positive integer degree. There is a universal
singular frame in which all divergences disappear. When computed
as iterated integrals, its coefficients are certain rational
numbers that appear in the local index formula of
Connes--Moscovici. When working with formal Laurent series over
$\Q$, the data of equisingular flat vector bundles define a
Tannakian category whose properties are reminiscent of a category
of mixed Tate motives.
\end{abstract}

\section{Introduction}

In this paper we show that the divergences of quantum field theory
are a highly structured phenomenon. More precisely, they provide
data that define an action of a specific {\em ``motivic Galois
group"} $U^*$ on the set of physical theories.

\smallskip

In particular, this exhibits the renormalization group  as
the action of a  one parameter subgroup
$\bG_a \,\subset  U^* $
of the above Galois group.

\smallskip

The work of Connes-Kreimer \cite{CK}, \cite{cknew} provided a
conceptual understanding of perturbative renormalization in terms
of the Birkhoff decomposition of loops in a pro-unipotent Lie
group $G$ determined by the physical theory, through the Hopf
algebra of Feynman graphs \cite{dhopf}, \cite{CK}.

\smallskip

This suggests the possibility of formulating the theory of
renormalization in the context of the Riemann--Hilbert
correspondence. The latter is a broad term encompassing, in
various forms and levels of generalization, equivalences between
geometric problems associated to differential systems with
singularities and representation theoretic data associated to the
monodromy.

\smallskip

In this paper we construct the Riemann--Hilbert correspondence
associated to perturbative renormalization, in the form of a
classification of flat equisingular bundles in terms of
representations of the ``motivic Galois group'' $U^*$.

\smallskip

More specifically, we start by considering the scattering formula
 \begin{equation} \g_- (z) = \lim_{t \to \infty}
e^{-t \left( \frac{\b}{z} + Z_0 \right)} \, e^{t Z_0} \,
\label{scattering}
 \end{equation}
proved in \cite{cknew}, which expresses the counterterms through
the residues of graphs.

\smallskip

We re-express this formula in terms of  the time ordered
exponential of physicists (also known as expansional in
mathematical terminology). The expression in expansional form can
be recognized as solution of a differential system. This
identifies a class of connections naturally associated to the
differential of the regularized quantum field theory viewed as a
function of the complexified dimension. The physics input that the
counterterms are independent of the additional choice of a unit of
mass translates, in geometric terms, into the notion of
equisingularity for these connections.

\smallskip

Thus, the geometric problem consists of the classification of {\em
``equisingular"} $G$-valued flat connections on the total space
$B$ of a principal $\bG_m$-bundle over an infinitesimal punctured
disk $\Delta^*$. An equisingular connection is a $\bG_m$-invariant
$G$-valued connection, singular on the fiber over zero, and
satisfying the following property: the equivalence class of the
singularity of the pullback of the connection by a section of the
principal $\bG_m$-bundle only depends on the value of the section
at the origin.

\smallskip

This classification problem
stems directly from the divergences of the physical theory at the
dimension $D$ where one would like to do physics\footnote{
We may assume $D=4$ (no strings attached).}.
The base $\Delta^*$ is the space of complexified dimensions around
$D$. The fibers of the principal $\bG_m$-bundle $B$
describe the arbitrariness in the normalization
of integration in complexified dimension $z\in \Delta^*$,
 in the commonly used
regularization procedure known as Dim-Reg (dimensional
regularization). The  $\bG_m$-action corresponds to the rescaling
$\hbar\, \partial /\partial \hbar$. The group $G$ is the
pro-unipotent Lie group whose Hopf algebra is the Hopf algebra of
Feynman graphs of \cite{dhopf}, \cite{CK}.

\smallskip

On the other side of our Riemann--Hilbert correspondence, the
representation theoretic setting equivalent to the classification
of equisingular flat connections is provided by representations
$U^* \to G^*$, where $U^*$ is a universal group, unambiguously
defined independently of the physical theory. The group $G^*$ is
the semi-direct product of $G$ by the action of the grading
$\t_t$, as in \cite{cknew}. We give an explicit
description of $U^*$ as the semi-direct product by its grading of
the graded pro-unipotent Lie group $U$ whose Lie algebra is the
free graded Lie algebra $$\cF(1,2,3,\cdots)_{\bullet}$$ generated
by elements $e_{-n}$ of degree $n$, $n>0$. 

\smallskip

Thus, there are three different levels at which Hopf algebra structures
enter the theory of perturbative renormalization. First, there is
Kreimer's Hopf algebra of rooted trees \cite{dhopf}, which is adapted
to the specific physical theory by decorations of the rooted
trees. There is then the Connes--Kreimer Hopf algebra of Feynman
graphs, which is dependent on the physical theory by construction, but
which does not require decorations. There is then the algebra
associated to the group $U^*$, which is universal with respect
to the set of physical theories.

\smallskip

We then construct a specific universal singular frame on principal
$U$-bundles over $B$. When using in this frame the dimensional
regularization technique of QFT, all divergences disappear and one
obtains a finite theory, which only depends upon the choice of a
local trivialization for the principal $\bG_m$-bundle $B$.

\smallskip

The coefficients of the universal singular frame, written out in
the expansional form, are the same that appear in the local index
formula of Connes--Moscovici \cite{cmindex}. In particular, they
are rational numbers. This means that we can view equisingular
flat connections on finite dimensional vector bundles as endowed
with arithmetic structure. We show that they can be organized into
a Tannakian category with a natural fiber functor to the category
of vector spaces, over any field of characteristic zero. The
Tannakian category obtained this way is equivalent to the category
of finite dimensional representations of the affine group scheme
$U^*$, which is uniquely determined by this property.

\smallskip

Closely related group schemes appear in
motivic
Galois theory and $U^*$
is for instance abstractly (but non-canonically)
isomorphic to
the motivic Galois group $ G_{\cM_T}(\fO) $ (\cite{dg}, \cite{sasha4}) of
 the scheme $S_4={\rm Spec}(\fO) $
of $4$-cyclotomic integers, $\fO=\,\Z[i][\frac{1}{2}]$.

\smallskip

The natural appearance of the ``motivic Galois group'' $U^*$ 
in the context of renormalization 
confirms a suggestion made by Cartier in \cite{Cart1},
that in the Connes--Kreimer theory of perturbative
renormalization one should find a hidden ``cosmic Galois group''
closely related
in structure to the Grothendieck--Teichm\"uller group. The question of
relations between the work of Connes--Kreimer, motivic Galois theory,
and deformation quantization was further emphasized by Kontsevich in
\cite{Kont}. At the level of the Hopf algebra of rooted trees,
relations between renormalization and motivic Galois theory were
also investigated by Goncharov in \cite{Gon02}.

\smallskip

The ``motivic Galois group'' $U$ acts on the set of dimensionless
coupling constants
of physical theories, through the map of the corresponding group $G$
to formal diffeomorphisms constructed in \cite{cknew}.

\smallskip

This also realizes the hope formulated
in \cite{CoGal1} of relating concretely the
renormalization group to a Galois group.
Here we are dealing with the Galois group
dictated by renormalization
and the renormalization group appears
as a canonical one parameter subgroup
$\bG_a \,\subset  U $.

\smallskip

These facts altogether indicate that the divergences
of Quantum Field Theory, far from just
being an unwanted nuisance, are a clear sign of the presence of
totally unexpected symmetries of geometric origin.
This shows, in particular, that one should
understand how the universal singular frame
``renormalizes" the geometry of space-time
using the Dim-Reg scheme and the universal
counterterms.

\bigskip

\section{Expansional form of the counterterms}\label{expan}

\bigskip

The following discussion will be quite general. We let $G$ be
a complex graded pro-unipotent Lie group, $\fg=\hbox{Lie} \ G$ its Lie
algebra, and $\t_t=e^{tY}$ the one parameter group of automorphisms
implementing the grading $Y$. We assume that the grading
$Y$ is integral and strictly positive.

\smallskip

We let $G^*$ be the semi-direct product
 \begin{equation}
 G^* = G \, \rtimes_{\t} \, \R  \label{Ren15}
 \end{equation}
 of $G$ by the action of the grading $\t_t$, hence
the Lie algebra of $G^*$ has an  additional generator $Z_0$, such
that
 \begin{equation}
 [Z_0 , X] = Y(X) \qquad \forall \, X \in \hbox{Lie} \ G \,. \label{Ren13}
 \end{equation}

\smallskip

We let $\Hc$ be the commutative Hopf algebra of coordinates on $G$.
For any unital algebra $A$ over $\C$, we let $G(A)$ be the group of points
of $G$ over $A$ \ie of homomorphisms
$$
\Hc \to   A\,,
$$
with the product coming from the coproduct of $\Hc$.

\smallskip

We identify the elements of the
Lie algebra $\fg=\hbox{Lie} \ G$
with linear forms $L$ on $\Hc$ such that
$$
L(X\,Y)=\, L(X)\,\ve(Y) +\, \ve(X)\, L(Y)\,,\quad \forall X\,,Y \in \Hc\,,
$$
where $\ve$ is the augmentation of $\Hc$, playing the role of the
unit in the dual algebra. More generally,
for any unital algebra $A$ over $\C$, one defines $\fg(A)$ as
the Lie algebra of linear
maps $\Hc \to  A$, fulfilling the above derivation rule.

\smallskip

In \cite{cknew} a complete characterization is given of those
$G$-valued loops $\g_\mu(z)$ satisfying the properties
\begin{equation}
  \g_{e^t \mu} (z) = \t_{t z} (\g_{\mu} (z)) \qquad \fl \, t \in
 \R \, ,  \label{Ren6}
 \end{equation}
 \begin{equation}
 \frac{\partial}{\partial \mu} \, \g_{\mu}^- (z) = 0 \,
 . \label{Ren3} \end{equation}
Here $\g_{\mu}^-$ is the negative part of the Birkhoff
decomposition
\begin{equation}
\g_\mu \, (z) = \g_{\mu}^- (z)^{-1} \, \g_{\mu}^+ (z) \qquad z \in
\partial\Delta \,
,\label{renorm15}
\end{equation}
where $\g_{\mu}^+$ and $\g_{\mu}^-$ extend to holomorphic maps on
$\Delta$ and $\P^1(\C)\smallsetminus \{ 0
\}$, respectively.

\smallskip

In this Birkhoff decomposition $\g_{\mu}^+$ provides the
renormalized values at $D$ and $\g_{\mu}^-$ provides the
counterterms for the renormalization procedure of
quantum field theory (\cf \cite{cknew}).
The properties \eqref{Ren6} and \eqref{Ren3} originate from
physical considerations, namely from the fact that the counterterms
are independent of the choice of the mass scale
parameter $\mu$ (\cf \cite{Collins} 7.1.4 p. 170).

\smallskip

We can regard the $\g_{\mu}$ as elements of $G(K)$, where we
let $K$ be the field $ \C(\{z\})$ of convergent Laurent series in $z$.

\smallskip

Given a $\fg=\hbox{Lie} \ G$-valued smooth function $\a(t)$
 where $t\in[a,b]\subset \R$ is a real parameter,
one defines the expansional (\cf \cite{Araki}), or time ordered
exponential,  by the equality
\begin{equation}\label{expansional}
 {\bf {\rm T}e^{\int_a^b\,\a(t)\,dt}}=\,1+\, \sum_1^\infty \int_{a\leq
s_1\leq \cdots\leq
s_n\leq b} \,\a(s_1)\cdots\,\a(s_n) \prod ds_j \, ,
 \end{equation}
where the  product comes from the coproduct in $\Hc$.

\smallskip

This defines an element of $G(\C)$, which is the value $A(b)$ at $b$
of the unique solution $A(t)$ with $A(a)=1$ at $t=a$ of the
differential equation
\begin{equation}\label{diffexp}
 dA(t)=\,A(t)\,\a(t)\,dt\,.
\end{equation}

\smallskip

The basic property of the expansional is the identity
\begin{equation}\label{expan1}
 {\bf {\rm T}e^{\int_a^c\,\a(t)\,dt}}=\,{\bf {\rm
T}e^{\int_a^b\,\a(t)\,dt}}\,
{\bf {\rm T}e^{\int_b^c\,\a(t)\,dt}}.
 \end{equation}
\medskip

\smallskip

With this notation, we can rewrite the scattering formula
\eqref{scattering} as follows.

 \begin{The}\label{genmu}  Let $\g_\mu(z)$ be a family of $G$-valued loops
 fulfilling \eqref{Ren6} and \eqref{Ren3}.
 Then there exists uniquely $\beta \in \fg$
 and a loop $\g_{\rm reg}(z)$ regular  at $z=0$
 such that
\begin{equation}\label{gmuexp}
  { \g_\mu(z) }=\,{\bf {\rm T}e^{-\frac{1}{z}\,
 \int^{-z \log\mu}_\infty\,\t_{-t}(\beta)\,dt}}\;
  \t_{z\log\mu}(\g_{\rm reg}(z))\,.
\end{equation}
Conversely, given any $\beta$ and any regular loop  $\g_{\rm reg}(z)$,
the expression \eqref{gmuexp}
gives a solution to  equations \eqref{Ren6} and \eqref{Ren3}.
  \end{The}

The Birkhoff
 decomposition of the loop $\g_\mu(z)$ of \eqref{gmuexp}
is given by
\begin{equation}\label{Birkhexp}
\begin{array}{rl}
 \g_\mu^+(z)= & {\bf {\rm T}e^{-\frac{1}{z}
  \,\int_0^{-z \log\mu}\,\t_{-t}(\beta)\,dt}}\;
 \t_{z\log\mu}(\g_{\rm reg}(z))\,, \\[3mm]
  \g_\mu^-(z) = & {\bf {\rm
 T}e^{-\frac{1}{z}\,\int_0^\infty\,\t_{-t}(\beta)\,dt}}\,.
\end{array}
\end{equation}

\bigskip

\section{Local equivalence of meromorphic connections}

\bigskip

We consider the local behavior, on an infinitesimal punctured disk
$\Delta^*$ centered at $z=0$, of solutions of $G$-differential
systems.

\smallskip

As above, we work with convergent Laurent series. Namely, we let
$K$ be the field $ \C(\{z\})$
of convergent Laurent series in $z$ and
$O \subset   K$ be
the subring of  series without a pole at $0$.
The field $ K$ is a differential field and we let
$\Omega^1$ be the $1$-forms on $ K$
with
$$
d\,: K \to \Omega^1
$$
the differential, $df = \frac{df}{dz}\,dz$.

\smallskip

A connection on the trivial principal $G$-bundle
$P= \Delta^*\times G$ is specified by the restriction of
the connection form to $\Delta^*\times 1$,
\ie by a
$\fg$-valued $1$-form $\omega$ on $\Delta^*$.
We let
$\Omega^1(\fg)$
denote $\fg$-valued $1$-forms on $\Delta^*$,
so that every element of $\Omega^1(\fg)$
is of the form $A \,dz$ with
$A \in \fg(K)$.

\smallskip

The operator
\begin{equation}\label{Doper}
 D: G( {K}) \to \Omega^1(\fg) \ \ \ \ Df=
f^{-1}\, df
\end{equation}
satisfies
\begin{equation}\label{prodrule}
 D(fh)=D h + h^{-1}\, Df\, h.
\end{equation}

\smallskip

We consider differential equations of the form
\begin{equation}\label{Gdiff}
 Df=\,\omega ,
\end{equation}
where $\omega\in \Omega^1(\fg)$ specifies the connection on the
trivial principal $G$-bundle.

\smallskip

\begin{Def}\label{equivABdef}
We say that two connections $\omega$ and $\omega'$ are equivalent iff
\begin{equation}\label{equivAB}
 \omega' = Dh + h^{-1} \omega h,
\end{equation}
for some $h \in  G(O)$.
\end{Def}

This simply identifies connections that differ by a change of local
frame, given by a $G$-valued map regular in $\Delta$.

\smallskip

By construction, the group $G$ is a projective limit of
linear algebraic groups $G_i$ whose Hopf algebras
 are finitely
generated graded Hopf subalgebras $\Hc_i \subset \Hc$.
Given $\omega\in\Omega^1(\fg)$,
its projections $p_i(\omega)\in\Omega^1(\fg_i)$
have a positive radius of convergence $\rho_i >0$.
Thus, for a choice of a base point $z_0\neq 0$ with $\vert
z_0\vert<\rho_i$, we obtain the monodromy in the form
\begin{equation}\label{Mono}
M=\, {\bf {\rm T}e^{\int_0^1\,c^*(\omega)}} ,
\end{equation}
where $c(t)$ is a simple closed path of winding number
one in the punctured disk of radius $\rho_i$, with endpoints
$c(0)=z_0=c(1)$.

\smallskip

When passing to the projective limit, one has to take care
of the change of base point, but the
triviality of the monodromy, $M=1$, is a well
defined condition. It ensures the existence
of solutions $f\in G(K)$ for equation \eqref{Gdiff}.

\smallskip

A solution $f$ of \eqref{Gdiff} defines a $G$-valued
loop. By our assumptions on $G$, any
$f\in G(K)$ has a unique
Birkhoff decomposition of the form
\begin{equation}\label{fBirkh}
f = (f^-)^{-1} f^+,
\end{equation}
with
$$
f^+ \in  G(O)\,,\quad f^-\in G(\scr{Q})
$$
where $O \subset  {K}$ is the subalgebra of regular functions
and $\scr{Q}=\,z^{-1}\,\C([z^{-1}])$. Since $\Qc$ is not unital one needs
to be more precise in defining $G(\Qc)$. Let
$\tilde{\Qc}=\,\C([z^{-1}])$ and $\ve_1$ its augmentation. Then
$G(\Qc)$ is the subgroup of $G(\tilde{\Qc})$ of homomorphisms
$\phi\,:\,\Hc \mapsto \tilde{\Qc}$ such that $\ve_1\circ
\phi=\,\ve$ where $\ve$ is the augmentation of $\Hc$.

\smallskip

\begin{Pro}  \label{class} Two connections $\omega_1$ and $\omega_2$
with trivial monodromy
are equivalent iff solutions $f_j$ of $ Df=\,\omega_j
$  have the same negative part in the
Birkhoff decomposition,
$$f^-_1=f^-_2\,.$$
\end{Pro}

\bigskip

\section{Classification  of equisingular flat
connections}\label{classif}

\bigskip

We now modify the geometric setting of the previous section, by
introducing a principal $\bG_m$-bundle
\begin{equation}\label{base}
\bG_m \to B\,\to \Delta\,,
\end{equation}
over the infinitesimal disk $\Delta$. We let
$$
 b\mapsto w(b) \quad \forall w\in \C^*
$$
denote the action of the multiplicative group $\bG_m=\C^*$.
We let $\pi\,:B\to \Delta$ be the
projection, with
$$V=\pi^{-1}(\{0\})\subset B$$
the fiber over $0\in \Delta$ and
$y_0 \in V$ a base point. We let
$B^*\subset B$ denote the complement of $V$.

\smallskip

We consider again a group $G$ as above, with grading $Y$.
We can then view the trivial principal $G$-bundle
$P=B\times G$ as equivariant with respect to $\bG_m$,
 using the action
\begin{equation}\label{Gdiff1}
 u(b,g)=\,(u(b),u^Y(g)) \quad \forall u\in \C^*\,,
\end{equation}
\medskip
where $u^Y$ makes sense, since the grading $Y$ is
integer valued.

\smallskip

\begin{Def}
Let $P^*=B^*\times G$ be the restriction to $B^*$
of the bundle $P$.
We say that a connection $\omega$ on $P^*$
is {\em equisingular} if it is $\bG_m$-invariant
and its restrictions to sections of the
principal bundle $B$ that agree at $0\in \Delta$
are all equivalent in the sense of Definition \ref{equivABdef}.
\end{Def}

\medskip

We consider again the operator $Df=
f^{-1}\, df$ as in \eqref{Doper} satisfying \eqref{prodrule}.
We have the following notion of equivalence for $G$-differential
systems on $B$.

\begin{Def}\label{equivP}
Two connections $\omega$ and $\omega'$ on $P^*$
are equivalent iff
$$ \omega' = Dh + h^{-1} \omega h, $$
for a $G$-valued $\bG_m$-invariant map $h$
regular in $B$.
\end{Def}

\smallskip

The main step towards the formulation of perturbative
renormalization as a Riemann-Hilbert correspondence is
given by the following correspondence
between flat equisingular $G$-connections and elements in
the Lie algebra $\fg$.

\smallskip

We begin by choosing a non-canonical regular section
$$
\sigma\,: \Delta \to B\,, \quad \text{ with } \quad
\sigma(0)= y_0 \,.
$$
We shall later show that the correspondence established below
is in fact independent of the choice of $\sigma$.
To lighten notations we use $\sigma$ as the local frame that
trivializes the bundle $B$, which we identify with $\Delta \times
\C^*$.

\bigskip
\begin{center}
\includegraphics[scale=0.35]{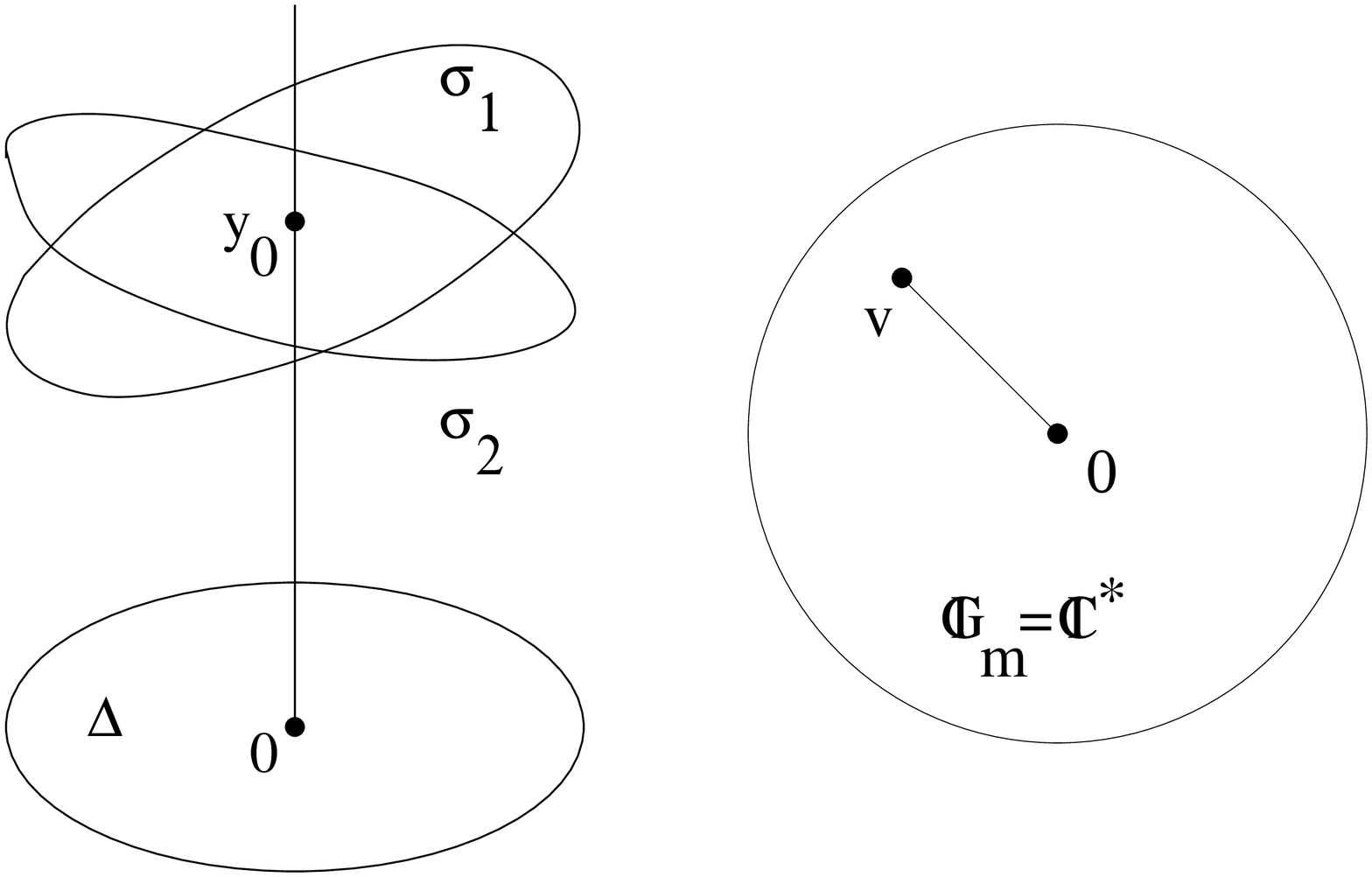}
\end{center}
\bigskip

\begin{The} \label{rh} Let $\omega$ be a flat
 equisingular $G$-connection.
There exists a unique element
 $\beta \in \fg$ of
the Lie algebra of $G$, such that $\omega$ is equivalent
to the flat
 equisingular connection $D\g$ associated to the section
\begin{equation}\label{solexpm}
 \g(z,v) =\,{\bf {\rm T}e^{-\frac{1}{z}\,
 \int^{v}_0\,u^Y(\beta)\,\frac{du}{u}}}\;
\in G\;
  \,,
\end{equation}
where the integral is performed on the straight path
$u=t v$, $t\in[0,1]$.
\end{The}

\proof
as above, we express a connection on $P^*$ in terms
of $\fg$-valued $1$-forms on $B^*$, and we
use the trivialization $\sigma$ to write it as
$$
\omega =\,A\,dz\,+\,B\,\frac{dv}{v},
$$
where both $A(z,v)$ and $B(z,v)$ are $\fg$-valued functions
 and $\frac{dv}{v}$ is the
fundamental $1$-form of the principal bundle $B$.

\smallskip

Let $\omega =\,A\,dz\,+\,B\,\frac{dv}{v}$ be
an invariant connection. One has
$$
\omega(z,u\,v)=\,u^Y(\omega(z,v))\,,
$$
which shows that $\omega$ is determined by its
restriction to  the section  $v=1$.
 One then has
$$\omega(z,u)=\,u^Y(a)\,dz\,+\,u^Y(b)\,\frac{du}{u} $$
for suitable elements $a,b\in \fg(K)$.

\smallskip

The flatness of the connection means that we have
\begin{equation}\label{flat}
\frac{db}{dz}-Y(a)+\,[a,b]=0.
\end{equation}

\smallskip

The positivity of the integral grading $Y$ shows that
the connection $\omega$ extends to a flat connection
on the product $\Delta^*\times \C$.
Moreover its restriction to $\Delta^*\times \{0\}$
is equal to $0$. This suffices to show that
the connection has trivial monodromy with respect to
both generators of $\pi_1(B^*)=\Z^2$.

\smallskip

One can then explicitely
write down a solution of
the differential system
\begin{equation}\label{partial}
D\g =\, \omega
\end{equation}
in the form
\begin{equation}\label{sol}
 \g(z,v) =\,{\bf {\rm T}e^{
 \int^{v}_0\,u^Y(b(z))\,\frac{du}{u}}}\;
  \,,
\end{equation}
where integration is performed on the straight path
$u=t v$, $t\in[0,1]$.

\smallskip

This gives a translation invariant loop $\g$,
\begin{equation}\label{gsol}
\g(z,u)=\,u^Y\g(z)
\end{equation}
fulfilling
\begin{equation}\label{abpartial}
\g(z)^{-1}d\g(z)=\,a\,dz\,,\qquad \g(z)^{-1}Y\g(z)=\,b\,.
\end{equation}

\smallskip

By hypothesis $\omega$ is equisingular,
hence the restrictions $\omega_s$ to the
lines $$\Delta_s^*=\,\{(z,e^{sz});\,z\in \Delta^* \}$$
are mutually equivalent. By proposition \ref{class}, using
the fact that the restriction of
$\g(z,u)=u^Y\g(z)$ to $\Delta_s^*$ is given
by $\g_s(z)=\,\t_{sz}\g(z)$, we obtain  that the
negative parts of the Birkhoff decomposition of
the loops $\g_s(z)$ are independent of the
parameter $s$.

\smallskip

Thus, by the results of Section \ref{expan},
there exists an element $\beta\in \fg$ and a regular
loop $\g_{\rm reg}(z)$, such that
\begin{equation}\label{solexp}
 \g(z,1) =\,{\bf {\rm T}e^{-\frac{1}{z}\,
 \int^{0}_\infty\,\t_{-t}(\beta)\,dt}}\;
 \g_{\rm reg}(z)\,.
\end{equation}

\smallskip

Thus, up to an equivalence given by the regular
loop $u^Y(\g_{\rm reg}(z))$, we can write
the solution in the form
\begin{equation}\label{solexp1}
 \g(z,u) =\,u^Y({\bf {\rm T}e^{-\frac{1}{z}\,
 \int^{0}_\infty\,\t_{-t}(\beta)\,dt}})\;
  \,,
\end{equation}
which only depends upon $\beta \in \fg$.
Since $u^Y$ is an automorphism one can
in fact rewrite \eqref{solexp1}
 as
\begin{equation}\label{solexp2}
  \g(z,v) =\,{\bf {\rm T}e^{-\frac{1}{z}\,
 \int^{v}_0\,u^Y(\beta)\,\frac{du}{u}}}\;
  \,,
\end{equation}
where the integral is performed on the straight path
$u=t v$, $t\in[0,1]$.

\smallskip

We then need to understand in what way the
class of the solution \eqref{solexp1}
depends upon $\beta \in \fg$.

\smallskip

An equivalence between two equisingular flat
connections generates a relation between
solutions of the form
$$
\g_2(z,u) =\,\g_1(z,u) \,h(z,u)
$$
with $h$ regular. Thus, the negative parts
of the Birkhoff decomposition of both
$$
\g_j(z,1) =\,{\bf {\rm T}e^{-\frac{1}{z}\,
 \int^{0}_\infty\,\t_{-t}(\beta_j)\,dt}}
$$
have to be the same, and this gives $\beta_1=\beta_2$.

\smallskip

Finally, we need to show that, for any $\beta \in \fg$,
the connection $\omega=\,D\g$ with $\g$
given by \eqref{solexpm}
is equisingular. The equivariance
follows from the invariance of the section $\g$.
Let then $v(z)\in\C^*$ be a regular
function with $v(0)=1$ and consider
the section $v(z)\sigma(z)$ instead
of $\sigma(z)$. The restriction of $\omega=\,D\g$
to this new section is now given by $D\g_v$, where
\begin{equation}\label{solexpm1}
 \g_v(z) =\,{\bf {\rm T}e^{-\frac{1}{z}\,
 \int^{v(z)}_0\,u^Y(\beta)\,\frac{du}{u}}}\;
\in G\;
  \,.
\end{equation}

\smallskip

We claim that the Birkhoff decomposition of $\g_v$
is given
by $\g_v(z) =\,\g_v^-(z)^{-1}\,\g_v^+(z)$ with
\begin{equation}\label{solexpm2}
\g_v^-(z)^{-1} = \, {\bf {\rm T}e^{-\frac{1}{z}\,
 \int^{1}_0\,u^Y(\beta)\,\frac{du}{u}}} \ \ \ \text{ and } \ \ \ 
\g_v^+(z) =\,{\bf {\rm T}e^{-\frac{1}{z}\,
 \int^{v(z)}_1\,u^Y(\beta)\,\frac{du}{u}}}\, .
\end{equation}
Indeed, the first term in \eqref{solexpm2} is
a regular function of $z^{-1}$ and
gives a polynomial in $z^{-1}$ when paired with
any element of $\Hc$. The second term is
a regular function of $z$, using the Taylor
expansion of $v(z)$ at $z=1$.

\endproof

\smallskip

By a similar argument, one gets the independence on the choice of
the section, as follows.

\begin{The} \label{ind} The above correspondence
between  flat
 equisingular $G$-connections
and elements
 $\beta \in \fg$ of
the Lie algebra of $G$ is independent of the choice of
the local regular section
$\sigma\,: \Delta \to B$, with $\sigma(0)= y_0$.

Given two choices  $\sigma_2=\,\alpha \,\sigma_1$
of local sections, the regular values
$\g_{{\rm reg}}(y_0)_j$ of solutions of  the differential
system above, in the corresponding singular frames, are related by
$$
\g_{{\rm reg}}(y_0)_2=\,e^{-s\,\beta}\g_{{\rm reg}}(y_0)_1
$$
where $$s= \,\left(\frac{d\alpha(z)}{dz}\right)_{z=0}\,.$$
\end{The}

\bigskip

It is this second statement that controls the ambigu\"ity
inherent to the renormalization group action in the physics setting,
where there is no preferred choice of
local regular section $\sigma$. In that context the principal bundle
$B$ over an infinitesimal disk of complexified dimensions around $D$
admits as fiber over $z\in \Delta$ the set of all possible
normalizations for the integration in complexified dimension $D-z$.
Moreover, the choice of the base
point in the fiber $V$ over $D$ corresponds to
the choice of the Planck constant, while the
choice of the section $\sigma$ (up to order one) corresponds
to the choice of a ``unit of mass".

\bigskip
\section{The universal singular frame}
\bigskip

We shall now reformulate the results of Section \ref{classif}
as a Riemann-Hilbert correspondence.
At the representation theoretic level, we
want to encode the data classifying equivalence
classes of equisingular flat connections
(Theorem \ref{rh}) by a homomorphism
\begin{equation}\label{repUG}
U^* \longrightarrow G^*
\end{equation}
from some universal group $U^*$ to $G^*$.

\smallskip

Viewed in this perspective, the group $U^*$ can be thought of as an
analog of the Ramis exponential torus in the
wild fundamental group that gives the local
Riemann--Hilbert correspondence in the context of differential Galois
theory (\cf \cite{Ramis}, \cite{vdp}). In fact, here the
equisingular flat connections have trivial monodromy and
one does not see the Stokes phenomenon, as we are only dealing with
perturbative renormalization. Thus, the group $U^*$ resembles most the
remaining part of the wild fundamental group, given by the
exponential torus, which appears in the formal local theory (\cf
\cite{Ramis}, \cite{vdp} and \S 3 of \cite{PuSi}). We will analyze
more closely the relation to the wild fundamental group in \cite{CM3}.

\smallskip

In \eqref{repUG} we need to get both $Z_0$ and $\beta$ in the range
at the Lie algebra level. Thus, working with Lie algebras, it
is natural to consider first the free Lie algebra
generated by $Z_0$ and $\beta$.
It is important, though, to keep track of the properties
one needs so that the formulae above make sense, such as
positivity  and integrality of the grading.

\smallskip

By these properties, we can
write $\beta$ as
an infinite formal sum
\begin{equation}\label{betasum}
\beta=\; \sum_1^\infty\;\beta_n\,,
\end{equation}
where, for each $n$, $\beta_n$ is homogeneous
of degree $n$ for the grading,
$$ Y(\beta_n)=n\, \beta_n. $$

\smallskip

Thus, assigning $\beta$ and the action of the grading on it
is the same as giving a collection of
homogeneous elements $\beta_n$ fulfilling
no restriction besides $Y(\beta_n)=n \beta_n$.
In particular, there is no condition on their
Lie brackets. Thus, these data are the same
as giving a homomorphism from the following
affine group scheme $U$ to $G$.

\smallskip

At the Lie algebra level, $U$ comes from the
free graded Lie algebra
$$ \cF(1,2,3,\cdots)_{\bullet} $$
generated by elements $e_{-n}$ of degree $n$, $n>0$.
At the Hopf algebra level, we therefore take the graded dual
of the enveloping algebra $\cU(\cF)$, so that
\begin{equation}\label{hopfu}
\Hc_u=\; \cU(\cF(1,2,3,\cdots)_{\bullet})^\vee .
\end{equation}

\smallskip

It is well known that, as an algebra, $\Hc_u$ is isomorphic
to the linear space of noncommutative polynomials
in variables $f_n$, $n\in \N_{>0}$ with the shuffle product.

\smallskip

We have obtained this way a pro-unipotent affine group scheme
$U$ which is graded in positive degree.

\smallskip

We can now reformulate the main theorem of Section \ref{classif}
as follows.

\begin{The} \label{rh1} Let $G$ be a positively graded
pro-unipotent Lie group.
There exists a canonical bijection between equivalence classes
of flat equisingular $G$-connections
and graded representations
$$\rho \, : U\,\to G $$
in $G$ of the group scheme $U$ defined above.
\end{The}

\smallskip

We can consider the
semi-direct product $U^*$ of $U$ by the grading as an affine
group scheme with a natural homomorphism
$U^*\,\to \bG_m$ to the multiplicative group.
The compatibility with the grading means that
$\rho$ extends to a homomorphism
$$\rho^* \, : U^*\,\to G^* .$$

\smallskip

Theorem \ref{rh1} shows that the group $U^*$ plays, in the formal
theory, a role analogous to the Ramis exponential torus of differential
Galois theory. The conceptual reason for considering the group $U^*$
rather than $U$ will become clear in the next section.

\smallskip

The equality
\begin{equation}\label{esum}
e=\; \sum_1^\infty\;e_{-n}
\end{equation}
defines an element of the Lie algebra of $U$.
Since $U$ is by construction a pro-unipotent
affine group scheme, we can lift $e$ to
a morphism of affine group schemes
\begin{equation}\label{rgmap}
{\bf{rg}}\; :\,\bG_a \,\to \,U\,,
\end{equation}
from the additive group $\bG_a$ to $U$.

\smallskip

The morphism \eqref{rgmap} represents the renormalization group in our
context. The corresponding ambigu\"ity is generated, as explained
above in Theorem \ref{ind}, by the absence of a canonical
trivialization for the $\bG_m$-bundle corresponding to integration
in complexified dimensions around $D$.

\smallskip

The formulae considered in the previous sections still make sense in
the universal case where $G^*=U^*$, hence we can define the
universal singular frame, by the equality
\begin{equation}\label{univ}
\g(z,v) =\,{\bf {\rm T}e^{-\frac{1}{z}\,
 \int^{v}_0\,u^Y(e)\,\frac{du}{u}}}\; \in U\; \,.
\end{equation}

\smallskip

This is easily computed in terms of iterated integrals and one obtains
the following expression.

\begin{Pro} \label{univsin} The
universal singular frame is given by
$$
\g(z,v) =\,\sum_{n \geq 0}\,\sum_{k_j>0}\,\frac{e(-k_1)e(-k_2)\cdots e(-k_n)}
{k_1\,(k_1+k_2)\cdots (k_1+k_2+\cdots +k_n)}\,v^{\sum k_j}\,z^{-n} .
$$
\end{Pro}

\smallskip

It is interesting to notice that exactly the same expression occurs in
the local index formula of \cite{cmindex}. The renormalization group
idea is also used in that context, in the case of higher poles in the
dimension spectrum.

\smallskip

Adopting this
universal singular frame in the dimensional regularization
technique has the effect of removing all divergences. One obtains a finite
theory, which depends only upon the choice of local
trivialization of the principal
$\bG_m$-bundle $B$, whose base $\Delta$ is the
space of complexified dimensions around $D$ and whose
fibers correspond to normalizations of the integral
in complex dimensions, as used by physicists
in the Dim-Reg scheme.

\bigskip
\section{The classifying affine group scheme
as a motivic Galois group}
\bigskip

In this section we construct a category of equivalence classes
of equisingular flat vector bundles. This allows us to reformulate
the Riemann--Hilbert correspondence in terms of finite dimensional
linear representations of $U^*$. The relation to the formulation given
in the previous section is given by passing to the finite dimensional
representations of the group $G^*$.
Since $G^*$ is an
affine group scheme, there are enough such representations,
and they are specified (\cf \cite{dg}) by assigning the following data.
\begin{itemize}
\item A graded vector space $E=\oplus_{n\in \Z} E_n$,
\item A graded representation $\pi$ of $G$ in $E$.
\end{itemize}

\smallskip

Notice that a graded representation of $G$ in $E$
can equivalently be described as a graded
representation of $\fg$ in $E$.  Moreover, since the
Lie algebra $\fg$ is positively graded, both representations
are compatible with the {\em weight} filtration
given by
\begin{equation} \label{weight}
W^{-n}(E)=\,\oplus_{m\geq n} E_m\,.
\end{equation}

\smallskip

At the group level, the corresponding representation in
the associated graded
$$
Gr^W_n=\, W^{-n}(E)/W^{-n-1}(E)\,.
$$
is the identity.

\smallskip

We now proceed to construct a Tannakian category of equivalence classes
of equisingular flat vector bundles, independent of the group $G$.

\begin{Def}\label{Wconn}
Let $(E,W)$ be  a filtered vector bundle
with a given trivialization of the associated graded
$Gr^W(E)$.
\begin{enumerate}
\item A $W$-connection on $E$ is a
connection $\nabla$ on $E$, which is compatible
with the filtration (\ie restricts to all $W^k(E)$)
and induces the trivial connection on the
associated graded $Gr^W(E)$.
\item Two $W$-connections on $E$ are $W$-equivalent iff there exists
an automorphism of $E$ preserving the filtration,
inducing the identity on $Gr^W(E)$, and conjugating
the connections.
\end{enumerate}
\end{Def}

\smallskip

We now define the category $\Ec$ of equisingular flat bundles.

\smallskip

Let $B$ be the principal $\bG_m$-bundle considered in Section
\ref{classif}.
The {\em objects} of $\Ec$ are the equivalence classes of pairs
$$\Theta=(E,\nabla) ,$$
where
\begin{itemize}
\item $E$ is a $\Z$-graded finite dimensional vector space.
\item $\nabla$ is an equisingular flat $W$-connection on $B^*$,
defined on the $\bG_m$-equivariant filtered vector bundle
$(\tilde{E},W)$ induced by $E$ with its weight filtration
\eqref{weight}.
\end{itemize}

\smallskip

By construction $\tilde{E}$ is the trivial bundle $B\times E$
endowed with the action of $\bG_m$ given by the grading.
The trivialization of the associated graded
$Gr^W(\tilde{E})$ is simply given by the identification
with the trivial bundle with fiber $E$.
The equisingularity of $\nabla$ here means that it is
$\bG_m$-invariant and that all restrictions to
sections $\sigma$ of $B$ with $\sigma(0)=y_0$
are $W$-equivalent.

\smallskip

We refer to such pairs $\Theta=(E,\nabla)$ as {\em flat equisingular
bundles}. We only retain the datum of the $W$-equivalence class of the
connection $\nabla$.

\smallskip

Given two flat equisingular bundles $\Theta$, $\Theta'$ we define the {\em
morphisms}
$$
T\in{\rm Hom}(\Theta, \Theta')
$$
in the category $\Ec$
as linear maps
$
T\,: \,E\to E'\,,
$
compatible with the grading,
fulfilling the condition that the following
$W$-connections $\nabla_j$, $j=1,2$, on $\tilde{E'}\oplus \tilde{E}$
are $W$-equivalent,
\begin{equation}\label{morphisms}
\nabla_1 = \left[ \begin{matrix}\nabla' &T\,\nabla-\,\nabla'\,T \cr
0 &\nabla \cr \end{matrix}
\right] \, \sim
\nabla_2 = \left[ \begin{matrix}\nabla' &0 \cr
0 &\nabla \cr \end{matrix}
\right] \,.
\end{equation}

\smallskip

In \eqref{morphisms}, $\nabla_1$
is obtained by conjugating $\nabla_2$ by the unipotent matrix
$$
\left[ \begin{matrix}1 &T \cr
0 &1 \cr \end{matrix}\right]\,.
$$
This shows that condition \eqref{morphisms}  is well defined, 
independently of the choice of representatives for
the connections $\nabla$ and $\nabla'$.

\smallskip

For  $\Theta=(E,\nabla)$, we set
$\omega(\Theta)=E$ and we view $\omega$
as a functor from the category of equisingular flat bundles
to the category of vector spaces. We then have the following result.

\smallskip

\begin{The} \label{tann}
Let $\Ec$ be the category of equisingular flat bundles defined above.
\begin{enumerate}
\item $\Ec$ is a Tannakian category.
\item The functor $\omega$ is a fiber functor.
\item $\Ec$ is equivalent to the category of finite dimensional
representations of $U^*$.
\end{enumerate}
\end{The}

\medskip

In all the above we worked over $\C$, with convergent Laurent
series. However, much of it can be rephrased with formal
Laurent series. Since the universal singular frame is
given in rational terms by proposition \ref{univsin}, the result
of Theorem \ref{tann}
 holds
over any field of characteristic zero and in
particular over $\Q$.

\smallskip

For each integer $n \in \Z$, we then define an object $\Q(n)$ in the
category $\Ec$ of equisingular flat bundles as the trivial bundle
given by a one-dimensional $\Q$-vector space placed in degree $n$,
endowed with the trivial connection on the associated bundle over $B$.

\smallskip

For any flat equisingular bundle $\Theta$ let
$$
\omega_n(\Theta)=\,{\rm Hom}(\Q(n), Gr^W_{-n}(\Theta))\,,
$$
and notice that $\omega = \oplus \,\omega_n$.

\bigskip

The group $U^*$ can be regarded as a motivic Galois group.
One has, for instance, the following identification  (\cite{sasha4},
\cite{dg}).

\begin{Pro} \label{motgal} There is a
(non-canonical)  isomorphism
\begin{equation}\label{MotU}
U^* \sim G_{\cM_T}(\fO) \,.
\end{equation}
of the affine group scheme $U^*$ with the motivic Galois group
$G_{\cM_T}(\fO) $
 of  the scheme $S_4$ of $4$-cyclotomic integers.
\end{Pro}

\smallskip

It is important here to stress the fact (\cf the ``mise en garde" of
\cite{dg}) that there is so far no ``canonical" choice
of a free basis in the Lie algebra of the above
motivic Galois group so that the above isomorphism
still requires making a large number of non-canonical
choices. In particular it is premature to assert that
the above category of equisingular flat bundles is
directly related to the category of $4$-cyclotomic
Tate motives. The isomorphism \eqref{MotU} does
not determine the scheme $S_4$ uniquely. In fact, a
similar isomorphism holds with $S_3$ the scheme
of 3-cyclotomic integers.

\smallskip
On the other hand, when considering the
category $\cM_T$ in relation to physics, inverting the prime $2$ is
relevant  to the definition of geometry
in terms of $K$-homology, which is at the center stage
in noncommutative geometry. We recall, in that respect, that it is
only after inverting the prime $2$
that (in sufficiently high dimension) a manifold structure on a simply
connected homotopy type is determined by the
$K$-homology fundamental class.

\smallskip

Moreover, passing from $\Q$ to a
field with a complex place, such
as the above cyclotomic fields $k$, allows for the
existence of non-trivial regulators for all
algebraic $K$-theory groups $K_{2n-1}(k)$.
It is noteworthy also that algebraic K-theory and
regulators already appeared in the context of
quantum field theory and NCG in \cite{coka}.
The appearance of multiple polylogarithms in
the coefficients of divergences in QFT,
discovered by Broadhurst and Kreimer (\cite{B}, \cite{BK}), as well as
recent considerations of Kreimer on analogies between residues of
quantum fields and variations of mixed Hodge--Tate structures
associated to polylogarithms (\cf \cite{Kr2}),
suggest the existence for the above
category of equisingular flat bundles of
suitable Hodge-Tate realizations given by
a specific choice of Quantum Field Theory.

\end{document}